\documentclass[12pt,a4paper]{article}
\pdfoutput=1
\usepackage[english]{babel}
\usepackage[utf8]{inputenc}
\usepackage[T1]{fontenc}
\usepackage{amsmath}
\usepackage{hyperref}
\usepackage{amssymb}
\usepackage{amsfonts}
\usepackage{amsthm}
\usepackage{multirow,tabularx}
\usepackage{graphicx}
\usepackage{setspace}
\usepackage{csquotes}

\usepackage{hyperref}
\usepackage{blkarray}
\usepackage{tikz}
\usepackage{lipsum}
\usepackage{setspace}
\usepackage{xspace}
\theoremstyle{plain}
\newtheorem{mydef}{Definition}
\newtheorem{theorem}{Theorem}[section]
\newtheorem{corollary}{Corollary}[theorem]
\newtheorem{lemma}[theorem]{Lemma}

\usepackage[reftex]{theoremref}
\usepackage{biblatex}
\bibliography{main.bib}

\title{Spectra and the exact number of  Automorphisms of Paley-type graphs of order a product of 'n' distinct primes}
\author{
 A. Sivaranjani\\
  \texttt{sivaranjani.a2021@vitstudent.ac.in}\\
  \text{Vellore Institute of Technology}\\
  \text{Chennai, 600 127}
  \and
  S. Radha\\
  \texttt{radha.s@vit.ac.in}\\
  \text{Vellore Institute of Technology}\\
  \text{Chennai, 600 127}
}
\date{August 2023}
\begin{document}
\maketitle

\begin{abstract}

   Paley graphs are Cayley graphs which are circulant and strongly regular. Paley-type graph of order a product of two distinct Pythagorean primes was introduced by Dr Angsuman Das. In this paper, we extend the study of Paley-type graphs to the order a product of ‘n’ distinct Pythagorean primes $p_1<p_2<p_3<...<p_n$.We have determined its adjacency spectra and the exact number of the automorphisms.
\end{abstract}
$\begin{aligned}
\end{aligned}$
\section[20pt]{Introduction}
Matrix theory and linear algebra were initially utilised to analyse networks using adjacency matrices. Algebraic approaches are very useful when dealing with regular and symmetric graphs. The study of the relationship pertaining to combinatorial features of graphs and the eigenvalues of matrices associated with the graph is known as Spectral Graph theory. \\ \\
Cayley graphs, named after mathematician Arthur Cayley, are a notion in the relationship between the concepts of groups and graphs. When a group G is represented as a Cayley graph, properties such as its size and the number of generators become considerably easier to investigate.
Cayley graphs are particularly important in the subject of combinatorics due to their ability to express the algebraic structure of groups. 
Cayley graphs are built to visualise cyclic groups $Z_n$, dihedral groups $D_n$, symmetric groups $S_n$, alternating groups $A_n$, direct and semidirect products and other finite group representations. Cayley graphs can help in understanding algebraic structures of groups in relation to representation theory. Refer  \cite{2} for more details on the characteristics of Cayley graphs and their many applications across a variety of fields.\\
\\ \\
Also, numerous distinct features of the Kronecker product of graphs have been examined (in different names, such as  direct product, cardinal product, conjunction, tensor product,  etc.). The study of structural findings is cited in  [\cite{6},\cite{4}, \cite{7}, \cite{8}, \cite{9}],their Hamiltonian features in directed\cite{11} and undirected graphs [\cite{10}, \cite{12}],their hyperbolicity  (see \cite{5}) stands out above all. A more in-depth structural understanding of this product would be beneficial, according to open problems in the field.\newline\\
Paley graph is a Cayley graph in which $V=\mathbb F_q$ a finite field on $q=p^{n}$ such that $q\equiv 1\mod{4}$ and the generating set S is the set of all elements of quadratic residue modulo q. Paley graphs are strongly regular, self-complementary, vertex-transitive, arc-transitive and edge-transitive. Paley graphs find applications in Coding theory, network design for the use to model and analyze network topologies and in the field of quantum computing in error detection and correction,etc.\newline \\

In this paper, we extend the study of Paley-type graphs introduced by  Dr.Angsuman Das \cite{1} to the order a  product of 'n' distinct Pythagorean primes and studied some basic graph theoretical properties with some known existing results.
We have also found out the adjacency spectrum and the exact number of automorphisms of the Paley-type graph.
We have also proved that the Paley-type graph $\Gamma_N$ is pancyclic for $p_1 > 5$.
\section{Definitions and Preliminaries} 
For the convenience of the readers, we provide some fundamental definitions, notations and some results in the form of lemmas and theorems without proof.
\subsection{Basic Definitions}
\begin{mydef}Let G be a finite group and let $S\subset G$ such that
    \begin{itemize}
        \item $e\notin S $( e is the identity of G)
        \item If $a\in S$ then $a^{-1}\in S$ and
        \item S generates G
    \end{itemize}
The Cayley graph $\Gamma$ is defined as $\Gamma=(V,E)$ where $V=G$  and \\ $E=\{ab | a,b \in V$ and $ a*b^{-1} \in S\}$.\\
\end{mydef}
\begin{mydef}
    \textbf{(Neighborhood of a node v in $\Gamma$)}\newline     Let $\Gamma$ be a graph.The neighbourhood of the node $v \in \Gamma$ written as $N(v)$, refers to the collection of nodes adjacent to the node v. The closed neighbourhood of v, is defined as $N[v]=N(v)\cup v$. \\
\end{mydef}

\begin{mydef}
\textbf{(Spectrum of Paley graph) \cite{13}}\newline  The eigenvalues of a given undirected Paley graph $\Gamma_q$ on q vertices where  $q \equiv 1\mod 4$  are $\frac{1}{2}(q-1)$ with multiplicity 1 and $\frac{1}{2}(-1\pm\sqrt{q})$ each with multiplicity ${\frac{1}{2}(q-1)}$.\\
\end{mydef}
\begin{mydef}\label{spec}
\textbf{(Spectrum of the Kronecker product of Adjacency \\ matrices of two graphs) \cite{14}} 
\newline Suppose $\Gamma_1$ and $\Gamma_2$ are two graphs with orders m and n respectively. Let  \textbf{A}  be the adjacency matrices of $\Gamma_1$ with eigenvalues $\lambda_1,\lambda_2,...,\lambda_m$ and \textbf{B}  be the adjacency matrix of $\Gamma_2$ with eigenvalues $\mu_1,\mu_2,...,\mu_n$. The eigenvalues of the Kronecker product of A and B  i.e., $\textbf{A}\otimes\textbf{B}$ are $\lambda_i\mu_j$ where $i=1,2,...,m$ and $j=1,2,...,n.$\\
\end{mydef}
\pagebreak
\begin{mydef}
    \textbf{\cite{3}}
    A graph $\Gamma$ is said to be \textbf{prime} with respect to the Kronecker product if it has more than one node and $\Gamma\cong \Gamma_1\times \Gamma_2$ implies that either $\Gamma_1$ or $\Gamma_2$ equals $K_1^{s}$ where $K_1^{s}$ denotes an isolated node with loop.\newline
The expression $\Gamma\cong \Gamma_1\times \Gamma_2 \times...\times \Gamma_k$ in which each $\Gamma_i$ prime is called the \textbf{prime factorisation} of $\Gamma$.\\
\end{mydef}
\begin{mydef}
\textbf{\cite{16}}A graph $\Gamma$ is called \textbf{R-thin} if no two nodes $a,b \in \Gamma_N$ have the
same open neighborhood i.e., $N_\Gamma(a)=N_\Gamma(b)$ implies a=b.\\
\end{mydef}

\subsection{Preliminaries}
Here, we set up some standard notations from the number theory  
\begin{itemize} 
\item ${\mathbb{Z}_N}$ denotes the set which contains all  integers modulo N 
\item ${\mathbb{Z}^*_N}$ denotes the set which contains all units in ${\mathbb{Z}_N}$ 
\item ${\mathcal{QR}}_N$ denotes the set of all quadratic residues which are also units in ${\mathbb{Z}_N}$ 
\item ${\mathcal{QNR}}_N$ denotes the set of all non-quadratic residues which are also  units in ${\mathbb{Z}_N}$ 
\item ${\mathcal{J}}^{+1}_N$ denotes the set of all values from ${\mathbb{Z}^*_N}$ with Jacobi symbol value +1 
\item ${\mathcal{J}}^{-1}_N$ denotes the set of all values from ${\mathbb{Z}^*_N}$ with Jacobi symbol  value -1 
\end{itemize} 

Despite the fact that this work is an expansion of \cite{1}, the same graph characteristics apply for n different primes. The following lemmas can be proved with the help of the results from elementary number theory.
\begin{lemma}
    If $p_1,p_2,...,p_n$  are 'n' distinct primes with $p_1\equiv p_2\equiv...\equiv p_n \equiv 1 (mod4)$, then -1 is a quadratic residue in $\mathbb{Z}_N$ 
    \end{lemma} 
\begin{theorem} \label{iso}
$\Gamma_N$ is isomorphic to the Kronecker product of Paley graphs $\Gamma_{p_1},\Gamma_{p_2},...,\Gamma_{p_n}$,where $p_i\equiv 1\mod 4$, i=1,2,...,n  i.e.,$\Gamma_N \cong \Gamma_{p_1} \times \Gamma_{p_2}\times...\times \Gamma_{p_n}$
\end{theorem}
\begin{lemma}{\label{diff}} 
\textbf{\cite{1}}Let $p$ be a Pythagorean prime and $z\in \mathbb{Z}_p$ .Then the number of ways in which $z$ can be represented as a difference of two quadratic residues in ${\mathbb{Z}^*_p}$ are
\begin{itemize}
    \item[1.] $\frac{p-1}{2}$ if $z\equiv(0\mod p)$ 
    \item[2.] $\frac{p-5}{4}$ if $z \in \mathcal{QR}_p$
    \item[3.] $\frac{p-1}{4}$ if $z \in \mathcal{QNR}_p$ 
\end{itemize}
\end{lemma}
\section{Paley-type graph and its symmetrical properties}
We first define the Paley-type graph on the product of 'n' distinct primes and we have listed some of their basic properties through lemmas and theorems. The proofs of the lemmas and theorems can be obtained using the same procedure given in \cite{1} 
\begin{mydef}
\textbf{ (Paley-type graph modulo N)}\newline For $N=p_1p_2...p_n$ , $p_{1}<p_2<...<p_n $ the Paley-type graph modulo $\Gamma_N$ is defined as $\Gamma_N=(V,E)$ where $V=\mathbb{Z}_N$ and $E=\{(a,b)\hspace{0.1cm}|\hspace{0.1cm} a-b\in \mathcal{QR}_N\}$ 
\end{mydef}
\textbf{Remark 1}:The graph $\Gamma_N$ is a Cayley graph (G,S) where the group $G=(\mathbb{Z}_N,+)$ and  the generating set $S=\mathcal{QR}_N$ \\ \\
\textbf{Remark 2}:   Throughout this paper we assume $p_i<p_{i+1}$, ${i=1,2,...,n}$

\begin{lemma}
    If $N=p_1p_2...p_n$, then the following holds:
\begin{itemize}
\item $\mathcal{QR}_N$ is a subgroup of ${\mathcal{J}}^{+1}_N$ and ${\mathcal{J}}^{+1}_N$ is a subgroup of ${\mathbb{Z}^*_N}$
\item $|{\mathbb{Z}^*_N}|=\phi(N)=(p_{1}-1)(p_{2}-1)...(p_{n}-1)$ 
\item $|\mathcal{QR}_N|=\frac{\phi(N)}{2^n}$ 
\item $|\mathcal{J}^{+1}_N|= 
|\mathcal{J}^{-1}_N|=\frac{\phi(N)}{2}$ 
\item $x\in \mathcal{QR}_N \iff 
x\in \mathcal{QR}_{p_{1}}\cap \mathcal{QR}_{p_{2}}\cap...\cap \mathcal{QR}_{p_{n}}$
    \end{itemize}
\end{lemma} 
\begin{theorem}
    The graph $\Gamma_N$ is Hamiltonian and hence connected.
    \end{theorem} 
    \begin{theorem}
    The graph $\Gamma_N$ is regular with degree $\frac{\phi(N)}{2^n}$ and hence Eulerian.
\end{theorem} 
\begin{theorem}
   The graph $\Gamma_N$ is both vertex and edge transitive.
    \end{theorem}
    \begin{theorem}
    Vertex and edge connectivity of the graph $\Gamma_N$ is given by $\frac{\phi(N)}{2^n}$.
    \end{theorem}
    \noindent
\textbf{Note:} $\Gamma_N$ is not self-complementary and not strongly regular.
\pagebreak
\begin{lemma} 
Let $N=p_1p_2\dots p_n$ where $p_1,p_2,\dots,p_n$ are 'n' distinct  
  primes. Then 
  \begin{itemize}
\item[1.] If $z\in{\mathcal{QR}_{N}}$ then $z$ can be represented as a difference of two quadratic residues i.e., $z={u^2}-{v^2}; u,v\in{{\mathbb{Z}}^*_N}$ in $$\frac{1}{4^n}\prod^{n}_{j=1} {(p_{j}-5)}$$ \newline
number of ways \\
\begin{proof}
    If $z \in \mathcal{QR}_N$ then $z \in \mathcal{QR}_{{p}_{j}}$. Hence, the proof follows by the Chinese remainder theorem and by (2) of Lemma \eqref{diff}.\\
    \end{proof} 
    
\item[2.]If $z \in {\mathcal{J}}^{+1}_N$ $\backslash$ $\mathcal{QR_{N}}$ and $z\in {{{\mathcal{QNR}_p}_j}_m},j_m \in \{1,2,...,n\} $; $r\leq n$; r is even, then the number of ways in which $z$ can be represented as difference of two quadratic residues given under two cases: \newline \\
Case (i): if $r=n$, then  $$\frac{1}{4^n}\prod^{n}_{j_m=1} {(p_{j_m}-1)}$$\\
Case(ii): if $r<n$, then
$$\frac{1}{4^n} \left[\prod^{r}_{m=1} ({p_{j_m}}-1) \prod^{n}_{k=r+1} {{{(p_j}_k}-5)}\right] $$ \\
\begin{proof} 
 Case (i): If $z \in {\mathcal{J}}^{+1}_N$ $\backslash$ $\mathcal{QR_{N}}$ and $z\in {{\mathcal{QNR}_p}_j}_m$ ; $r=n$ then by Chinese remainder theorem and by (3) of Lemma \eqref{diff},  the proof follows.\\ \\ \\

Case(ii): If $z \in {\mathcal{J}}^{+1}_N$ $\backslash$ $\mathcal{QR_{N}}$ and $z\in {{\mathcal{QNR}_p}_j}_m$; $r<n$ and even, then  by the Chinese remainder theorem and by (3) and (2) of Lemma\eqref{diff} the proof follows.\\
\end{proof}
\item[3.] If $z\in {\mathcal{J}}^{-1}_N$ and $z\notin {{\mathcal{QR}_p}_j}_k, j_k \in \{1,2,...,n\} $; $r\leq n;$ r is odd,  then the number of ways in which $z$ can be represented as difference of two quadratic residues is given under two cases:\newline \\
Case (i): if $r=n$, then  $$\frac{1}{4^n}\prod^{n}_{j_k=1} {(p_{j_k}-1)}$$
\newline \\
{Case (ii):} if $r<n$, then $$\frac{1}{4^n} \left[ \prod^{r}_{k=1} {{(p_j}_k}-5) \prod^{n}_{m=r+1} {{{(p_j}_m}-1)}\right] $$  \\
\begin{proof}
     Case (i): If $z\in {\mathcal{J}}^{-1}_N$  and $z\notin {{\mathcal{QR}_p}_j}_k$ ,$r=n$ then by Chinese remainder theorem and by (3) of Lemma \eqref{diff}, we complete the proof.\newline \\
Case(ii):  If $z\in {\mathcal{J}}^{-1}_N$ and $z \notin {{\mathcal{QR}_p}_j}_k$ ; $r<n$ and odd then by applying (2) of Lemma \eqref{diff} for ${{p_j}_k}$ and (3) of lemma \eqref{diff} for ${{p_j}_m}$ and by using Chinese remainder theorem we get the proof.\newline
\end{proof}
\item[4.]If $z(\neq 0) \in {\mathbb{Z}_N} \backslash {\mathbb{Z}}^*_N$ i.e., z is a  non-unit non-zero in ${\mathbb{Z}_N}$ then\newline \\
Case (i): For any $p_1,p_2,\dots,p_{n}$ if $z \equiv 0 \mod \prod^{n-1}_{m=1} {{p_j}_m} $ and $ z \in {{\mathcal{QR}_p}_j}_k;  k\in \{1,2,...n\};\hspace{0.2cm} j_m\neq j_k$ then 'z' can be represented as a difference of two quadratic residues is 
     $${\left[ \prod^{n-1}_{m=1} \frac{({{p_j}_m}-1)} {2^{n-1}}\right]} 
    {\left[\frac{({{p_j}_k}-5)}{4}\right]}$$ \newline
     number of ways 
\\ \\
Case(ii): If $z \equiv 0 \mod \prod^{r}_{m=1} {{p_j}_m}$ and $z \in {{\mathcal{QNR}_p}_j}_k;j_m\neq j_k$ where  $r<n;\hspace{0.2cm}  j_m,j_k \in \{1,2,...,n\}$ then $'z'$ can be represented as a difference of two quadratic residues is
     $${\left[ \prod^{n}_{j=1} \frac{(p_{j}-1)}{2^{r} .\hspace{0.1cm} 4^{n-r}} \right]}$$
     number of ways
\begin{proof}
Case(i): As $z \equiv 0 \mod \prod^{n-1}_{m=1} {{p_j}_m} $ and $ z \in {{\mathcal{QR}_p}_j}_k$ then by applying Chinese remainder theorem and by (1),(2) of Lemma \eqref{diff} we complete the proof. 
\\ \\
Case (ii): As $z \equiv 0 \mod \prod^{r}_{m=1} {{p_j}_m}$ and $z \in {{\mathcal{QNR}_p}_j}_k$ then by applying (1),(3) of lemma \eqref{diff} and by Chinese remainder theorem we end the proof. 
    \end{proof}
    \end{itemize}
    \end{lemma}
\section{Adjacency spectra of Paley-Type Graphs}
In this section, using the properties and results of the Kronecker product we have obtained the adjacency spectra of $\Gamma_N$.\newline 
\begin{theorem}
    Let $\Gamma_N$ be a Paley-type graph where $N=p_1p_2...p_n$.Then there are $3^n$ distinct eigenvalues in the adjacency spectra.
\end{theorem}
\begin{proof}
    Let us consider a Paley graph on two distinct primes $p_1,p_2$ say $\Gamma_{p_{1}p_{2}}$. Then the  adjacency spectrum of $\Gamma_{p_{1}p_{2}}$ can be obtained using the definition \eqref{spec} and is given in the following table.
   \\ \\ \\
\begin{tabular}{|m{4cm}|m{4cm}|m{3cm}|}\hline
\textbf{Eigenvalues}&\textbf{No of distinct eigenvalues}&\textbf{Multiplicities}\\ 
\hline \vspace{0.5cm} 
$\frac{(p_{1}-1)(p_{2}-1)}{2^2}$ \vspace{0.5cm}& \vspace{0.5cm} ${2^0}$ &\vspace{0.5cm} $1$ \vspace{0.5cm}\\
\hline \vspace{0.5cm}
$\frac{(\pm \sqrt p_{1}-1)(p_{2}-1)}{2^2}$ & \vspace{0.5cm} ${2^1}$ & \vspace{0.5cm} $\frac{(p_1-1)}{2}$\\  
\vspace{0.5cm} $\frac{(\pm \sqrt p_{2}-1)(p_{1}-1)}{2^2}$ & \vspace{0.5cm} ${2^1}$ &\vspace{0.5cm} $\frac{(p_2-1)}{2}$\\ 
\hline \vspace{0.5cm}
\vspace{0.5cm} $\frac{(\pm \sqrt p_1-1)(\sqrt p_2-1)}{2^2}$ & \vspace{0.5cm} ${2^2}$ &\vspace{0.5cm} $\frac{(p_1-1)(p_2-1)}{2}$\\ 
\vspace{0.5cm} $\frac{(\pm \sqrt p_1-1)(-\sqrt p_2-1)}{2^2}$ & \vspace{0.5cm} ${2^2}$ &\vspace{0.5cm} $\frac{(p_1-1)(p_2-1)}{2}$ \\& &\\ 
\hline 
\end{tabular}
\\ \\ \\ 
Therefore the number of distinct eigenvalues of  $\Gamma_{p_{1}p_{2}}$
\begin{equation*}
\begin{aligned}
&= 2_{C_0}\times 2^0+2_{C_1}\times 2^1+2_{C_2}\times 2^2\\
&
=(1+2)^2\\
&
=3^2\\
\end{aligned}
\end{equation*}
Extending the above to $N=p_1p_2...p_n$
we have summarised the spectrum of a Paley-type graph $\Gamma_N$ below.\\ \\ \\
\begin{tabular}{|m{6cm}|m{4cm}|m{3cm}|}\hline
\textbf{Eigenvalues}&\textbf{No of distinct eigenvalues}&\textbf{Multiplicities}\\ 
\hline \vspace{0.5cm} 
$\frac{(p_1-1)(p_2-1)...(p_n-1)}{2^n}$ \vspace{0.5cm}& \vspace{0.5cm} ${2^0}\vspace{0.5cm}$ &\vspace{0.5cm} $1$ \vspace{0.5cm}\\
\hline \vspace{0.5cm}
$\frac{(\pm \sqrt p_1-1)(p_2-1)...(p_n-1)}{2^n}$ & \vspace{0.5cm} ${2^1}$ & \vspace{0.5cm} $\frac{(p_1-1)}{2}$\\  
\vspace{0.5cm} $\frac{(p_1-1)(\pm \sqrt p_2-1)...(p_n-1)}{2^n}$ & \vspace{0.5cm} ${2^1}$ &\vspace{0.5cm} $\frac{(p_2-1)}{2}$\\  
\ldots & & \\
\ldots & & \\
\ldots & & \\
\vspace{0.5cm} $\frac{(p_1-1)( p_2-1)...(\pm \sqrt p_n-1)}{2^n}$ \vspace{0.5cm} &${2^1}$&$\frac{(p_n-1)}{2}$\\ \hline
\vspace{0.5cm}$\frac{(\pm \sqrt p_1-1)(\pm \sqrt p_2-1)...(p_n-1)}{2^n}$ &\vspace{0.5cm} ${2^2}$ & \vspace{0.5cm}$\frac{(p_1-1)}{2} \frac{(p_2-1)}{2}$ \\ 
\vspace{0.5cm} $\frac{(p_1-1)(\pm \sqrt p_2-1)(\pm \sqrt p_3-1)...(p_n-1)}{2^n}$ & \vspace{0.5cm}${2^2}$ & \vspace{0.5cm}$\frac{(p_2-1)}{2} \frac{(p_3-1)}{2}$ \\ 
\ldots & & \\ 
\ldots & & \\ 
\ldots  & & \\
$\frac{(p_1-1)(p_2-1)...(\pm \sqrt p_{n-1}-1)(\pm \sqrt p_n-1)}{2^n}$ & ${2^2}$ & $\frac{(p_{n-1}-1)}{2} \frac{(p_n-1)}{2}$\\
 ...&...&... \\
 ...&...&...\\ 
 ...&...&...\\ & &
\\$\frac{(\pm \sqrt p_1-1)(\pm \sqrt p_2-1)...(\pm \sqrt p_{n-1}-1)(\pm \sqrt p_n-1)}{2^n}$  &  ${2^n}$ &  $\frac{(p_{1}-1)(p_{2}-1)...(p_{n}-1)}{2^n} $\\& &\\  \hline 
\end{tabular}
\newline\\ \\
\\
Hence the number of distinct eigenvalues
\begin{equation*}
\begin{aligned}
&= n_{C_0}\times 2^0+n_{C_1}\times 2^1+n_{C_2}\times 2^2+...+n_{C_n}\times 2^n\\
&
=(1+2)^n\\
&
=3^n\\
\end{aligned}
\end{equation*}
\end{proof}
    \begin{corollary}
    The least eigenvalue of the Paley-type graph $\Gamma_N$ is  \\  $\frac{(-\sqrt{p}_1-1)(p_2-1)...(p_n-1)}{2^n}$ 
\end{corollary}
\section{Order of the Automorphism group of Paley-type Graph}
\begin{lemma}\label{rthinlemma}
    The Paley-type graph $\Gamma_N$ on $N=p_1 p_2...p_n$ vertices is R-thin.
\begin{proof}
 To prove this, let us assume  $N_{\Gamma_N}(a)=N_{\Gamma_N}(b)$ and prove $a=b$.\\
 Now
     $$N_{\Gamma_N}(a)=\{q+a| x^2\equiv q\mod N, x\in \mathbb{Z}_N,  q\in \mathcal{QR}_N\}$$
    $$N_{\Gamma_N}(b)=\{q+b| x^2\equiv q\mod N, x\in \mathbb{Z}_N,  q\in \mathcal{QR}_N\}$$
Then $$q+a=q+b$$ Implies $a=b$. Thus the graph $\Gamma_N$ is R-thin.
\end{proof}
\end{lemma}
\begin{lemma}\label{prime}
    If $\Gamma_N$ is prime.Then the Paley-type graph $\Gamma_N$ has  prime \\ factorisation $\Gamma_N \cong \Gamma_{p_1}\times\Gamma_{p_2}\times...\times\Gamma_{p_n}$
\end{lemma}
\begin{proof}
    By the definition 6, we can say that each $\Gamma_{p_i}$ can be written as \\ ${K_1}^s \times \Gamma_{p_i} \cong \Gamma_{p_i}$. Therefore each $\Gamma_{p_i}$ is prime. Hence $\Gamma_{p_1}\times\Gamma_{p_2}\times...\times\Gamma_{p_n}$ is a prime factorisation of $\Gamma_N$  .
\end{proof}
\begin{lemma}\label{aut}
    \cite{3} Suppose $\Gamma$ be a connected, non-bipartite and R-thin graph and let $\phi$ be an automorphism defined on it which has a prime factorisation. $$\Gamma=\Gamma_1\times \Gamma_2\times...\times \Gamma_k$$ Then there exists a permutation 
    $\tau$ of \{1,2,...,k\} together with isomorphisms \newline $\phi:\Gamma_{\tau(i)} \rightarrow \Gamma_i$ such that
$$\phi(x_1,x_2,...,x_k)=(\phi_1(x_{\tau(1)}),\phi_2(x_{\tau(2)}),...,\phi_n(x_{\tau(k)}))$$
Thus, the automorphism group of $\Gamma$ is generated by the automorphisms of the prime factors and transpositions of isomorphic factors. Hence, Aut( $\Gamma$) is isomorphic to the automorphism group of the disjoint union of the prime factors of $\Gamma$.\\
\end{lemma} 
\noindent
We use the above lemmas and propositions to prove the following theorem. \\
\begin{theorem}
    Let $Aut(\Gamma_N)$ denote the automorphism group of $\Gamma_N$. Then the number of automorphisms is $|Aut(\Gamma_N)|=\frac{N\phi(N)}{2^n}$
\end{theorem}
\begin{proof}
The Paley-type graph $\Gamma_N$ is connected and non-bipartite which is \\ obvious from the adjacency spectrum of $\Gamma_N$. The graph $\Gamma_N$ is R-thin by lemma \ref{rthinlemma} \hspace{0.1cm}Also each $\Gamma_{p_i}$ is prime and $\Gamma_N$ has a prime factorisation $\Gamma_{p_1}\times\Gamma_{p_2}\times...\times\Gamma_{p_n}$ by lemma \ref{prime}. Hence the graph $\Gamma_N$ is a connected,  non-bipartite R-thin graph and has a prime factorisation $\Gamma_{p_1}\times\Gamma_{p_2}\times...\times\Gamma_{p_n}$.\newline \\
Then, by lemma \ref{aut}
\begin{equation}\label{equ1}
    \begin{aligned}
|Aut({\Gamma_N})|= |Aut( {\Gamma_{p_1}})|   \times |Aut ({\Gamma_{p_2}})| \times \hdots \times |Aut({\Gamma_{p_n}})|  
\end{aligned}
\end{equation}
We know that the paley graph of order $q=p^e$ has  $|Aut(\Gamma_q)|=\frac{q(q-1)e}{2}$ \newline
For a Paley-type graph $q=p_i$ the order is $|Aut(\Gamma_{p_i})|=\frac{p_i({p_i}-1)}{2}$ \newline \\ 
    By \eqref{equ1}
   \begin{equation*}
\begin{aligned}
i.e.,|Aut(\Gamma_N)|&
=\frac{p_1(p_1-1)}{2}\times\frac{p_2(p_2-1)}{2}\times\hdots\times\frac{p_n(p_n-1)}{2}\\
&
=\frac{p_1\times p_2 \times \hdots \times p_n \times (p_1-1)\times(p_2-1)\times \hdots \times (p_n-1)}{2^n}
\\&
=\frac{N\phi(N)}{2^n}
  \end{aligned}
\end{equation*}
\end{proof}
\pagebreak
\section{Pancyclicity of Paley-type graphs}
\begin{lemma} \label{pan}
    The Kronecker product of two pancyclic graphs is pancyclic. 
\end{lemma}
\begin{proof}
    Let $G_1$ and $G_2$ be two pancyclic graphs of orders $n_1$ and $n_2$ respectively. Since $G_1$ and $G_2$ are pancyclic they have cycles of length k for $3\leq k \leq n_1$ and $3\leq k \leq n_2$. Consider k-cycles $x_1 x_2...x_k x_1 \in G_1$ and $y_1 y_2...y_k y_1 \in G_2$ . \\
    By the definition of the Kronecker product ,$(x_1,y_1) (x_2,y_2)...(x_k,y_k)(x_1,y_1)$ is a k-cycle in $G_1 \times G_2$. This is true for all k. Hence $G_1 \times G_2$ is pancyclic. 
\end{proof}
\begin{theorem} \cite{17}
    The Paley-type graph $\Gamma_N$ on $N=p_1p_2...p_n$, $p_1 > 5$ is pancyclic.
\end{theorem}
\begin{proof}
     $\Gamma_N$ is the Kronecker product of Paley graphs  $\Gamma_{p_1}\Gamma_{p_2}...\Gamma_{p_n}$ by the theorem \eqref{iso} . By the above lemma \eqref{pan} $\Gamma_N$ is pancyclic.
\end{proof}
\section*{Statements and Declarations}
The authors declare that no funds, grants or other support were received during the preparation of this manuscript.
\section*{Conclusion}
In this article, we have extended the study of Paley-type graphs and found out its adjacency spectrum, the exact number of automorphisms and the pancyclicity. Metric, local metric and partition metric dimensions of Paley-type graphs are still unresolved.
\printbibliography
\end{document}